\documentclass[a4paper,twocolumn]{autart-accepted}

\usepackage{calrsfs,times}

\usepackage{graphicx,psfrag,enumerate}
 \usepackage{pstool}
\usepackage{amssymb}
\usepackage{amsmath}
\usepackage{calrsfs}
\usepackage{color, epsfig, epstopdf, wrapfig}
\usepackage{paralist}
\usepackage{float,cuted}

\newtheorem{definition}{Definition}
\newtheorem{theorem}{Theorem}

\newtheorem{lemma}{Lemma} 
 
\newtheorem{assumption}{Assumption}
\newtheorem{proposition}{Proposition}

\usepackage{tikz-cd}
\usetikzlibrary{babel}

\usepackage{amsmath}
\usepackage{accents}
\newlength{\dhatheight}

\usepackage{color}
\usepackage[normalem]{ulem}
\usepackage{soul} 
\soulregister\cite7
\soulregister\ref7
\soulregister\pageref7

\begin{document}
\begin{frontmatter}

\title{A Robust Control Approach to Asymptotic Optimality of the Heavy Ball
  Method for Optimization of Quadratic Functions}

\thanks[footnoteinfo]{This work was 
    supported by the Australian Research Council under the Discovery Projects funding scheme (Project numbers DP190102158, DP200102945 and DP210102454). Part
    of this work was carried out during the first author's visit to the
    Australian National University.} 

\author[First]{V.~Ugrinovskii}\ead{v.ugrinovskii@gmail.com} \quad
\author[Second]{I.~R.~Petersen}\ead{i.r.petersen@gmail.com} \quad
\author[Second]{I.~Shames}\ead{iman.shames@anu.edu.au}

\address[First]{School of Engineering and IT, University of New South Wales
  Canberra, Canberra, ACT, 2600, Australia}
\address[Second]{The CIICADA Lab, School of Engineering, The Australian
        National University, Canberra, ACT 
        2601, Australia} 
       
\begin{abstract}
Among first order optimization methods, Polyak's heavy
ball method has long been known to guarantee the asymptotic rate of
convergence matching Nesterov's lower bound for functions
defined in an infinite-dimensional space. In this paper, we use results on
the robust gain margin of linear uncertain feedback control systems to show that the
heavy ball method is provably worst-case asymptotically optimal when
applied to quadratic functions in a finite dimensional space.  

\end{abstract}

\end{frontmatter}

\section{Introduction}
First order methods for solving optimization problems 
\begin{equation}
  \label{eq:66}
  \min_{x} f(x)
\end{equation}
iteratively approximate the minimum point $x^*$ of $f$ using linear combinations
of the previous iterates and the gradients of $f$ computed at those
previous iterates. Such methods find widespread applications in machine
learning and its applications to control. A fundamental problem regarding
such methods is to characterize their rate of convergence for a given
class of objective functions $f$~\cite{dAST-2021,Nesterov-2018,DT-2022}.

For two times differentiable unimodal functions
$f:\mathbf{R}^n\to \mathbf{R}$ whose  
Hessian at the stationary point $\Delta=\nabla^2 f(x^*)$ satisfies
\begin{equation}
  \label{eq:52}
mI\le\Delta\le LI,
\end{equation}
Polyak's heavy ball method~\cite{Polyak-1964a,Polyak-1987}, 
\begin{equation}
  \label{eq:5}
  x_{t+1}=x_t-\alpha\nabla f(x_t)+\beta(x_t-x_{t-1}), \quad t=0,1,\ldots,
\end{equation}
guarantees that
\begin{equation}
   r_{\{x_t\}}= \limsup_{t\to \infty} \|x_t-x^*\|^{1/t}\le \rho^*, \quad \rho^*\triangleq \frac{\sqrt{L}-\sqrt{m}}{\sqrt{L}+\sqrt{m}},
\label{eq:4}
\end{equation}
when it is initiated sufficiently close to $x^*$. Here $\nabla f(\cdot)$
denotes the gradient of $f$. The quantity $r_{\{x_t\}}$ is
the root-convergence factor of the sequence 
$\{x_t\}$~\cite{OR-2000}; it characterizes the asymptotic rate of
convergence of $\{x_t\}$ to $x^*$. When the set of
functions $f$ is restricted to include only quadratic functions
satisfying~(\ref{eq:52}) (we denote this class of functions
$\mathcal{Q}_{m,L}^n$), the method~(\ref{eq:5}) converges
globally~\cite{Polyak-1987}, and the inequality~(\ref{eq:4}) is 
tight in the worst-case:    
\begin{equation}
  \label{eq:6}
  \sup_{f\in\mathcal{Q}_{m,L}^n} r_{\{x_t\}}=\rho^*.
\end{equation}
Polyak commented~\cite[p.74]{Polyak-1987},  that for large-scale problems with
quadratic functions in $\mathcal{Q}_{m,L}^n$, 
where the dimension of the vector $x$ is greater than the number of
iterations $T$ required to reach $x^*$ with a sufficient accuracy, 
one cannot expect any first order method to converge at a rate faster than a
geometric sequence with the heavy ball method's ratio $\rho^*$. 

The same ratio $\rho^*$ appears in the (often misquoted)
infinite dimensional Nesterov's
lower bound. Nesterov showed~\cite[Theorem
2.1.13]{Nesterov-2018} that among $m$-strongly convex 
continuously differentiable functions with $L$-Lipschitz gradient, defined
in an infinite dimensional Hilbert space, there exists
a `bad' quadratic function $f$ for which 
\begin{equation}
  \label{eq:36}
  \|x_t-x^*\|\ge (\rho^*)^t\|x_0-x^*\|
  \quad \forall t=0,1,\ldots,
\end{equation}
for any first order optimization method. As a result, any such method
applied to this quadratic function will produce a sequence of iterates whose
root-convergence factor is bounded from below:  
\begin{equation}
  r_{\{x_t\}}=\limsup_{t\to \infty} \|x_t-x^*\|^{1/t}\ge \rho^*.
\label{eq:12}
\end{equation}
For functions in the finite dimensional space $\mathbf{R}^n$, lower bounds
similar to~(\ref{eq:36}) hold only over 
the first $T$ steps, and $T$ is linked to the dimension of
the space $\mathbf{R}^n$~\cite{dAST-2021,DT-2022}. These results have led
the community to
believe that in optimization of quadratic functions $\mathbf{R}^n\to 
\mathbf{R}$, no other first order method can provide a better  
root-convergence factor than the heavy ball method. However, this
conclusion cannot be drawn formally from the existing finite-dimensional
results, since they do not hold as $T\to\infty$ while the dimension $n$ of
the search space remains fixed. 

In this paper, we 
consider the set $\mathcal{Q}_{m,L}^n$ of finite-dimensional quadratic functions
$f:\mathbf{R}^n\to \mathbf{R}$ with the Hessian $\Delta$ satisfying
condition~(\ref{eq:52})
and fixed-parameter first-order methods
\begin{eqnarray}
x_{t+1}&=&x_t
+\sum_{j=0}^{k-1}\beta_j(x_{t-j}-x_{t-j-1})-\sum_{j=0}^{l}\alpha_j\nabla
           f(y_{t-j}), \nonumber
\\
y_t&=&\sum_{\nu=0}^{k-l}\gamma_\nu x_{t-\nu}, \quad t=0,1,\ldots.   \label{eq:50}
\end{eqnarray}
The number of past iterates $k$ and the number of gradient evaluations
$l\le k$ used at every step can be arbitrary, but they do not change with
time. The coefficients 
$\alpha_j$, $\beta_j$ and $\gamma_\nu$ are scalar constants. 

Our main result shows that for any converging 
method~(\ref{eq:50}), the worst-case
root-convergence factor within the 
class $\mathcal{Q}_{m,L}^n$ is
bounded from below by the same value $\rho^*$ which appears on the
right-hand side of equation~(\ref{eq:4}); see inequality~(\ref{eq:79.1}) in
Theorem~\ref{specrad}.
This conclusion applies to any
optimization algorithm~(\ref{eq:50}) including
the fixed-step gradient descent method~\cite{Nesterov-2018,Polyak-1987},
Polyak's heavy-ball method~\cite{Polyak-1964a,Polyak-1987}, the triple momentum
method~\cite{VFL-2017}, and Nesterov's fixed parameter accelerated
method~\cite{Nesterov-2018}. From this result
and equation~(\ref{eq:6}) it immediately follows that the heavy ball method
is \emph{worst-case optimal} among methods~(\ref{eq:50}) in the sense that for
quadratic 
functions in $\mathcal{Q}_{m,L}^n$ it guarantees the best worst-case
asymptotic convergence rate. 

Our approach uses robust control theory. Specifically, our derivation of
the worst-case lower bound on the root-convergence factor of the
method~(\ref{eq:50}) employs the results on the robust
gain margin of feedback control systems~\cite{KT-1985}.

\section{Optimal lower bound on the convergence rate of first order methods
  applied to quadratic functions}\label{degstab}

In this section, we derive the optimal lower bound on the 
root-convergence factor of the method~(\ref{eq:50}) applied to quadratic functions
$f$ of the class $\mathcal{Q}_{m.L}^n$. First we recall some basic
definitions which formalize the notion of the 
asymptotic rate of convergence of an iterative process. 

\begin{definition}[\cite{OR-2000}]\label{Def.R-convergence}
  Let $\{x_t\}$ be a sequence that converges to $x^*$. Then the number
  \begin{equation}
    \label{eq:28}
    r_{\{x_t\}}=\limsup_{t\to \infty} \|x_t-x^*\|^{1/t}
  \end{equation}
is the \emph{root-convergence factor}, or \emph{$R$-factor} of $\{x_t\}$.
If $\mathcal{I}$ is an iterative process with limit point
$x^*$, and $\mathcal{C}(\mathcal{I},x^*)$ is the set of all sequences
generated by 
$\mathcal{I}$ which converge to $x^*$, then 
\begin{equation}
  \label{eq:29}
  r_{\mathcal{I}}=\sup\{r_{\{x_t\}}: \{x_t\}\in \mathcal{C}(\mathcal{I},x^*)\}
\end{equation}
is the $R$-factor of $\mathcal{I}$ at $x^*$.
\end{definition}

Iterative processes of the form~(\ref{eq:50}) can be
written in the form of a nonlinear dynamic
system of Lur{\' e} type~\cite{HC-2008,Khalil},
  \begin{eqnarray}
    \label{eq:51}
&&X_{t+1}=AX_t+BU_t, \\
&& Y_t=CX_t,\quad U_t=-\phi(Y_t), \nonumber
  \end{eqnarray}
whose state, output and nonlinearity are respectively 
\begin{eqnarray*}
  && X_t= \left[
 \begin{array}{ccc}
 x_{t-k}^T & \ldots & x_t^T
 \end{array}
\right]^T, ~~
Y_t=
\left[
  \begin{array}{ccc}
   (C_l X_t)^T & \ldots & (C_0X_t)^T
  \end{array}
\right]^T, \\ 
&&\phi(Y)= 
    \left[\begin{array}{ccc}
    (\nabla f(C_lX))^T & \ldots &  (\nabla f(C_0X))^T 
          \end{array}
    \right]^T. \nonumber
  \end{eqnarray*}
The matrices $A\in \mathbf{R}^{(k+1)n\times (k+1)n}$, $B\in
  \mathbf{R}^{(k+1)n\times (l+1)n}$, $C\in \mathbf{R}^{(l+1)n\times (k+1)n}$ and $C_j\in \mathbf{R}^{n\times (k+1)n}$,
      $j=0,\ldots,l$, are defined as
\begin{eqnarray*}
&&A=A_0\otimes I_n, \quad B=B_0\otimes I_n, \quad  C=
   \left[
   \begin{array}{ccc}
     C_l^T& \ldots & C_0^T
   \end{array}
   \right]^T, \nonumber \\
&&A_0=   \left[
   \begin{array}{c|cccc}
\mathbf{0} & \multicolumn{4}{c}{I_{k}} \\ \hline 
-\beta_{k-1} & ~\beta_{k-1}-\beta_{k-2}~ &
  ~\beta_{k-2}-\beta_{k-3} & \ldots & 1+\beta_0
   \end{array}
   \right], \nonumber \\   
&&B_0=
   \left[
   \begin{array}{cccc}
\multicolumn{4}{c}{\mathbf{0}} \\ \hline
\alpha_l~ &~ \alpha_{l-1} & \ldots &~\alpha_0 
   \end{array}
   \right], \quad C_j=C_j^0\otimes I_n, \nonumber \\
&& C_j^0= 
   \underbrace{\left[\begin{array}{ccc}
0 & \ldots & 0
   \end{array}\right.}_{l-j} 
\begin{array}{ccc}
\gamma_{k-l} & \ldots & \gamma_0 
   \end{array}\underbrace{\left.\begin{array}{ccc}
0 & \ldots & 0
   \end{array}\right]}_j.
\label{eq:53}
  \end{eqnarray*}
In the above equations, the symbol $\otimes$ denotes the Kronecker product,
$I_n$ is the $n\times n$ identity matrix, and $\mathbf{0}$ denotes a zero
matrix or vector.
The dimensions of the matrices $A_0$, $B_0$ and $C_j^0$ are $(k+1)\times (k+1)$,
$(k+1)\times (l+1)$, and $1\times (k+1)$, respectively. 

\begin{assumption}\label{A1}
It holds that
\begin{equation}
  \label{eq:8}
\sum_{j=0}^l\alpha_j\neq 0, \quad \sum_{\nu=0}^{k-l}\gamma_\nu=1.
\end{equation}
\end{assumption}

Assumption~\ref{A1} holds for many practical methods. For instance,
in Nesterov's
constant-step accelerated gradient method~\cite[constant
step scheme III, p.94]{Nesterov-2018}, $l=0$, $k=1$, $\alpha_0=\frac{1}{L}\neq 0$, 
$\beta_0=\frac{1-\sqrt{m/L}}{1+\sqrt{m/L}}$, and $\gamma_0=(1+\beta_0)$,
$\gamma_1=-\beta_0$. Thus,  $\gamma_0+\gamma_1=1$. 

\begin{lemma}\label{L1}
Suppose $f(x)$ has a unique stationary point $x^*$. Under
Assumption~\ref{A1}, $x^*$  is the unique fixed point of the method
$\mathcal{I}$ in equation~(\ref{eq:50}).  
\end{lemma}

\emph{Proof: }
Substituting the $x_t=x_{t-1}=\ldots =x_{t-k}=x^*$ into the left-hand side
of equation~(\ref{eq:50}) and using the second identity~(\ref{eq:8}) and
the fact that $\nabla f(x^*)=0$ yields
\begin{eqnarray*}
 x_{t+1}&=&x^*- \sum_{j=0}^l\alpha_j\nabla
            f\left(\sum_{\nu=0}^{k-l}\gamma_\nu x^*\right) \\
&=& x^*-\left(\sum_{j=0}^l\alpha_j\right)\nabla
    f\left(\left(\sum_{\nu=0}^{k-l}\gamma_\nu\right)x^*\right) = x^*.
\end{eqnarray*}
This shows that $x^*$ is a fixed point of the mapping~(\ref{eq:50}).

To show uniqueness, suppose that $\check x$ is a
fixed point of the method $\mathcal{I}$. Then 
\[
\left(\sum_{j=0}^l\alpha_j\right)\nabla
    f\left(\left(\sum_{\nu=0}^{k-l}\gamma_\nu\right)\check x\right)=0.
\]  
According to~(\ref{eq:8}) this implies that $\nabla f(\check x)=0$,
i.e., $\check x$ is a
stationary point of $f$. Since $f$ is assumed to have a unique stationary
point, the fixed point of the method $\mathcal{I}$ must be unique.    
\hfill$\Box$

It follows from Lemma~\ref{L1} that when $f$ is unimodal,
$X^*=[\underbrace{1~\ldots~1}_{k+1}]^T\otimes x^*$ is the unique
equilibrium of the system~(\ref{eq:51}).  

For a quadratic function $f(x)=\frac{1}{2}(x-x^*)^T\Delta (x-x^*)+f_0$,
using the change of variable $\bar X_t= X_t-X^*$ and the condition
$\sum_{\nu=0}^{k-l}\gamma_\nu=1$ from Assumption~\ref{A1}, the
system~(\ref{eq:51}) can be written as the linear system
\begin{eqnarray}
    \label{eq:51lin}
&& \bar X_{t+1}=\bar A\bar X_t \quad \mbox{where }
  \label{eq:3}
\bar A\triangleq A_0\otimes I_n-(B_0\otimes \Delta)C.
  \end{eqnarray}
The system~(\ref{eq:51lin}) is stable if and only if the
spectral radius of the matrix $\bar A$, denoted $\rho(\bar A)$, satisfies
$\rho(\bar A)<1$. Moreover, $\rho(\bar A)$ characterizes the degree of
stability of the system~(\ref{eq:51lin})~\cite{AM90}. That is, when the
system~(\ref{eq:51lin}) is stable, the states $X_t$ approach $X^*$ at least
as fast as $\rho(\bar A)^t$. In terms of the root convergence factor, this
observation reads 
\begin{equation}
  \label{eq:13}
  r_{\mathcal{I}}=\rho(\bar A)
\end{equation}
for any $f\in\mathcal{Q}_{m,L}^n$ for which the matrix $\bar A$ is stable. 

When the function $f$ is not known in advance and it is known only that it
belongs to the set $\mathcal{Q}_{m,L}^n$, the 
system~(\ref{eq:51lin}) is uncertain. When it is stable for all $f\in
\mathcal{Q}_{m,L}^n$, its degree of stability in the face of this uncertainty
is characterized by $\sup_{mI\le \Delta \le LI} \rho(\bar A)$, and so the
states $X_t$ approach $X^*$ at least as fast as $(\sup_{mI\le \Delta \le
  LI} \rho(\bar A))^t$. That is, the worst-case root convergence rate of
the method~(\ref{eq:50}) is given by  
\begin{equation}
  \label{eq:14}
  \sup_{f\in \mathcal{Q}_{m,L}^n}r_{\mathcal{I}}= \sup_{mI\le \Delta \le LI} \rho(\bar A).
\end{equation}

We are now in a position to present the main result of the paper which
establishes the lower bound on the quantity
in~(\ref{eq:14}). Note that $\sup_{mI\le \Delta \le LI}\rho(\bar A)$ can be
greater than or equal to 1, however the expression on the left-hand side
of~(\ref{eq:14}) is well-defined only when the sequences $\{x_t\}$ and $\{X_t\}$
asymptotically converge. This is the standing assumption of this result.    

\begin{theorem}
  \label{specrad}
For any $k\ge 1$ and $0\le l\le k$ and any collection of parameters
$\boldsymbol{\alpha}=(\alpha_0,\ldots,\alpha_{l})$, 
$\boldsymbol{\beta}=(\beta_0,\ldots,\beta_{k-1})$ and 
$\boldsymbol{\gamma}=(\gamma_0,\ldots,\gamma_{k-l})$ which satisfies
Assumption~\ref{A1} and for which the algorithm~(\ref{eq:50})
asymptotically converges for all $f\in \mathcal{Q}_{m,L}^n$, it holds that 
\begin{equation}
  \label{eq:79.1}
  \sup_{f\in \mathcal{Q}_{m,L}^n}r_{\mathcal{I}}= \sup_{\Delta: mI\le
    \Delta\le LI} \rho(\bar A)\ge \rho^*.
\end{equation}
\end{theorem}

The proof of this theorem employs the result of~\cite{KT-1985} regarding
the robust gain margin of linear feedback control systems. It is presented
in Lemma~\ref{KT-lemma} given below. This lemma is a special case of a more
general theoretical development in~\cite{KT-1985}. To make our presentation
self-contained, we include a direct proof of Lemma~\ref{KT-lemma} in
the Appendix.

\begin{figure}[t]
\begin{center}
\psfragfig[width=0.9\columnwidth]{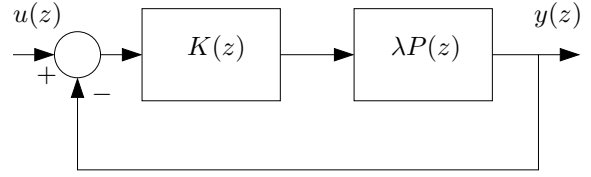}{ 
  \psfrag{C(z)}{\hspace{-1ex}$K(z)$}
  \psfrag{u}{$u(z)$}
  \psfrag{y}{$y(z)$}
  \psfrag{+}{$+$}
  \psfrag{-}{$-$}
  \psfrag{P(z)}{\hspace{-1ex}$\lambda P(z)$}}
\end{center}
  \caption{A linear uncertain control system with a robust stabilizing
    compensator.} 
  \label{fig:lincomp}
\end{figure}

\begin{lemma}
  \label{KT-lemma}
Consider the uncertain linear feedback system shown in Fig.~\ref{fig:lincomp}
consisting of the plant $P(z)=\frac{1}{z-1}$, an uncertain constant gain
$\lambda$
and a compensator $K(z)$. Let
$\rho\in (0,1)$ be a constant. A proper real rational compensator $K(z)$
that places 
all poles of this system in the interior of the disk $|z| < \rho$ for all
$\lambda \in [m,L]$ exists if and only if 
  \begin{equation}
    \label{eq:15}
    \rho> \rho^*,
  \end{equation}
where $\rho^*$ is the constant defined in~(\ref{eq:4}), i.e.,
$\rho^*\triangleq \frac{\sqrt{L}-\sqrt{m}}{\sqrt{L}+\sqrt{m}}$. 
\end{lemma}

{\emph{Proof of Theorem~\ref{specrad}: }
Since the algorithm~(\ref{eq:50}) is assumed to
converge asymptotically for all $f\in\mathcal{Q}_{m,L}^n$, then $\rho(\bar
A)<1$ for every symmetric matrix $\Delta$ which satisfies~(\ref{eq:52}). 
We now show that $\sup_{mI\le\Delta \le LI}\rho(\bar A)< 1$. Indeed, 
the set $mI\le \Delta\le LI$ is compact in the finite dimensional metric
space of symmetric matrices equipped with the norm induced by the Euclidean
norm in $\mathbf{R}^n$. This conclusion follows from the fact that this set is a
closed bounded set. Furthermore, $\rho(\bar A)$ depends continuously on
$\Delta$; see~(\ref{eq:3}). Thus, by the Weierstrass extreme value theorem,
there exists a symmetric matrix $\Delta^*$ within that set which attains the
supremum. However, for this matrix $\Delta^*$ the conditions of the
theorem state that the spectral radius of the corresponding
matrix~(\ref{eq:3}) $\bar A^*=A_0\otimes I_n-(B_0\otimes \Delta^*)C$ is
less than 1. Thus we conclude 
that  
\begin{equation}
  \label{eq:22}
\sup_{mI\le\Delta \le LI}\rho(\bar A)=\max_{mI\le\Delta \le LI}\rho(\bar A)
=\rho(\bar A^*)< 1.  
\end{equation}

Now consider an arbitrary function $f\in\mathcal{Q}_{m,L}^n$. The
corresponding matrix $\Delta$ is symmetric, therefore there exists an
orthogonal matrix $T$ such that $\Delta=T^T\Lambda T$, where $\Lambda$ is
the diagonal matrix whose diagonal consists of the eigenvalues of
$\Delta$. Using the matrix $T$, let us change coordinates in~(\ref{eq:51lin}),
\begin{eqnarray*}
  \widetilde X_t=\widetilde T \bar X_t \triangleq (I_{k+1}\otimes T)\bar X_t.
\end{eqnarray*}
In the new coordinates, the system~(\ref{eq:51lin}) becomes
\begin{eqnarray}
  \label{eq:56}
 &&\widetilde X_{t+1}=\widetilde T \bar A\widetilde T^T\widetilde X_t=\widetilde A\widetilde X_t,
\end{eqnarray}
where $\widetilde A=\widetilde T \bar A\widetilde T^T$.
The state transformation does not affect the spectral radius of a
matrix. Therefore, $\rho(\bar A)=\rho(\widetilde A)$.
Since the matrix $\widetilde T=I_{k+1}\otimes T$ is orthogonal, 
$\widetilde A =
A_0\otimes I-N \otimes \Lambda. 
$ 
Here $N$ is a $(k+1)\times(k+1)$ matrix of the form
\begin{eqnarray}
  \label{eq:58}
 N=
  \left[
  \begin{array}{cccc}
\multicolumn{4}{c}{\mathbf{0}} \\ \hline
    n_k & n_{k-1} &  \ldots  & n_0 
  \end{array}
  \right],
\end{eqnarray}
where $\mathbf{0}$ is the $k\times (k+1)$ zero matrix.
The elements of the last row of $N$ are the coefficients of the product
polynomial 
$(\sum_{j=0}^l\alpha_jz^j)(\sum_{\nu=0}^{k-l}\gamma_\nu z^\nu)$,   
\begin{equation}
  \label{eq:59}
  n_j=\sum_{\nu=0}^{j} \alpha_{\nu}\gamma_{j-\nu}, \quad j=0,\ldots,k.
\end{equation}
Here we use the standard convention that if the index extends
beyond the length of a vector, the corresponding element of the extended
vector is taken to be 0.

Let $\lambda_1,\ldots, \lambda_n$ be the eigenvalues of $\Delta$, and so
$\Lambda=\mathrm{diag}[\lambda_1,\ldots,\lambda_n]$. It is easy to see by
permuting the columns and rows of the matrix $A_0\otimes I-N\otimes
\Lambda$ that the eigenvalues of the matrix $\widetilde A$ are the
same as those of the block diagonal matrix comprised of the companion matrices
\begin{eqnarray*}
  \label{eq:60}
  g_i&=&A_0-\lambda_iN \nonumber \\
  &=&{\scriptsize
  \left[
  \begin{array}{cccc}
    0 & 1 & \ldots  & 0 \\
\vdots & \vdots & \ddots & \vdots \\ 
    0 & 0 & \ldots  & 1 \\
    -\beta_{k-1}-n_k\lambda_i & \beta_{k-1}-\beta_{k-2}-n_{k-1}\lambda_i &  \ldots  &1+\beta_0-n_0 \lambda_i
  \end{array}
  \right]}, \qquad \\
&& \hspace{5cm} i=1, \ldots,n.
\end{eqnarray*}
Therefore, we arrive at the conclusion that 
\begin{equation}
\label{eq:62}
\bar\rho\triangleq 
\sup_{mI\le \Delta\le LI}\rho(\bar A)=\max_i\sup_{\lambda_i\in [m,L]}\rho(g_i)
=\sup_{\lambda_i\in [m,L]}\rho(g_i).
\end{equation}
In the last identity, $\max_i$ is dropped since all eigenvalues $\lambda_i$
lie in the same interval $[m,L]$, and therefore $\bar\rho=\sup_{\lambda_i\in
  [m,L]}\rho(g_i)$ and does not depend on $i$.

Next we observe that $\rho(g_i)$ represents the radius of the
smallest disk in the complex plane which contains all roots of the
characteristic equation of the matrix $g_i$,
\begin{eqnarray}
  \label{eq:67}
(z-1)\left(z^k-\sum_{j=0}^{k-1}\beta_jz^{k-j-1}\right)+\lambda_i\sum_{j=0}^kn_jz^{k-j}=0.\quad
\end{eqnarray}
This equation can be written as
$
1+\lambda_iP(z)K(z)=0,
$
where $P(z)$ and $K(z)$ are given by
\begin{eqnarray*}
&&P(z)=\frac{1}{z-1}, \quad
K(z)=\frac{N(z)}{D(z)}, \nonumber \\
&&N(z)\triangleq\sum_{j=0}^kn_jz^{k-j}, \quad D(z)\triangleq
   z^k-\sum_{j=0}^{k-1}\beta_jz^{k-j-1}.
\nonumber
\end{eqnarray*}
Thus, we conclude that $\bar\rho=\sup_{\lambda_i\in [m,L]}\rho(g_i)$ is the
radius of the smallest disk that contains the poles of the SISO feedback
control  
system in Fig.~\ref{fig:lincomp} consisting of the uncertain plant $\lambda_i
P(z)$ and the compensator $K(z)$. According
to~(\ref{eq:22}), this radius $\bar\rho$ is in the interval $(0,1)$. Hence
one can select $\rho\in(\bar\rho,1)$ such that 
all poles of this system 
lie in the interior of the disk $|z|<\rho$ for any
$\lambda_i\in [m,L]$. 
From Lemma~\ref{KT-lemma},
a proper compensator $K(z)$ which ensures such pole placement for the family of
  plants $\lambda_iP(z)$ exists  
if and only if $\rho> \rho^*$. Taking infimum with respect to
$\rho\in(\bar\rho,1)$ proves (\ref{eq:79.1}).  
\hfill$\Box$
 
\section{Conclusions} 
This paper uses results on the robust gain margin of linear uncertain systems to
establish a theoretical lower bound on the rate of asymptotic 
convergence of iterative first order optimization methods applied to
quadratic functions. 
One conclusion that follows from our results is that among multistep
methods of the form~(\ref{eq:50}), Polyak's heavy ball
method~\cite{Polyak-1964a,Polyak-1987} in fact guarantees the best 
worst-case rate of asymptotic root convergence for the class of quadratic
functions. Also, our result confirms that Nesterov's bound indeed holds
asymptotically for functions defined in finite-dimensional space when one
considers the worst case over  quadratic functions. Consequently it also holds
for any class of strongly convex twice differentiable functions that are
$m$-convex and have $L$-Lipschitz gradient in the vicinity of $x^*$,
containing $\mathcal{Q}_{m,L}^n$ as a subset.  

\newcommand{\noopsort}[1]{} \newcommand{\printfirst}[2]{#1}
  \newcommand{\singleletter}[1]{#1} \newcommand{\switchargs}[2]{#2#1}

\section*{Appendix: Proof of Lemma~\ref{KT-lemma}}
Let $N(z)$, $D(z)$ be the numerator and denominator of the compensator
$K(z)$. Introduce the sensitivity function
\begin{eqnarray}
  \label{eq:16}
  S(z)&=&(1+\frac{m+L}{2}P(z)K(z))^{-1} .
\end{eqnarray}
Note that the pole of the plant $P(z)$, $z=1$, lies outside the disk
$|z|<\rho $ and is a zero of $S(z)$. Also, $z=\infty$, the
unique zero of $P(z)$, is a zero of $1-S(z)$ since the compensator
$K(z)$ is proper. Thus, we conclude that
\begin{equation}
  \label{eq:18}
S(1)=0, \quad S(\infty)=1.  
\end{equation}

The following proposition adapts Lemma~2.3 from~\cite{KT-1985} to the
problem setting of this paper. 

\begin{proposition}\label{KT-1985-L2.3}
  The closed loop system in Fig.~\ref{fig:lincomp} has all its poles in the
  disk $|z|<\rho$ for all $\lambda\in [m,L]$ if and only if 
$S(z)$ in~(\ref{eq:16}) is real rational and analytic in the region $\tilde
{\mathcal{H}_\rho}\triangleq \{|z|\ge\rho\}\cup\{\infty\}$ and
\begin{eqnarray}
  \label{eq:17}
  S(z)\not\in \mathcal{G}\triangleq \left(-\infty,\frac{2m}{m-L}\right]\cup
  \left[\frac{2L}{L-m},+\infty\right)\nonumber \\ \forall z\in \tilde
{\mathcal{H}_\rho}.
\end{eqnarray}
\end{proposition}

The proof of this proposition follows, \emph{mutatis mutandis}, the proof of
Lemma~2.3 in~\cite{KT-1985}. Based on this proposition, we conclude that
the poles of the closed loop system under 
consideration can be placed in the interior of the open disk $|z|<\rho$ for all
$\lambda\in [m,L]$ if and only if there exist polynomials $D(z)$,
$N(z)$ with real coefficients, such that the degree of
$N$ is less than or equal to the degree of $D$  and 
the function $S(z)$ in~(\ref{eq:16}) is analytic in $\tilde{\mathcal{H}}_\rho$, maps
$\tilde{\mathcal{H}}_\rho$ into the complement of the set
  $\mathcal{G}$, denoted $\mathcal{G}^c$ (i.e.,
$\mathcal{G}^c=\mathbf{C}\backslash\mathcal{G}$ where $\mathbf{C}$ is the
set of complex numbers), and satisfies~(\ref{eq:18}). 
Following~\cite{KT-1985}, we observe that the existence of such polynomials
and the function $S(z)$ is essentially the Nevanlinna-Pick interpolation
problem.  

The classical Nevanlinna-Pick problem~\cite{BGR-2013} is concerned with the
following. Given the points $\zeta_1, \zeta_2, \ldots \zeta_l$ in the
interior of the unit disk $\bar{\mathcal{D}}=\{z: |z|\le 1\}$ and the array
of values $b_1,b_2,\ldots, b_l$, the problem is to find a
rational function $s(z)$ with no poles in $\mathcal{D}=\{z: |z|< 1\}$, for
which $\sup_{z\in \mathcal{D}}|s(z)|<1$ and such that $s(\zeta_i)=b_i$,
$i=1,\ldots, l$.  
The equivalence between this classical formulation and the interpolation
problem stated above follows from the following diagram 
\[
\begin{tikzcd}
\tilde{\mathcal{H}}_\rho \arrow[r,"S"] \arrow[d,"\varphi"] &
  \mathcal{G}^c \arrow[d,"\theta"] \\
\bar{\mathcal{D}} \arrow[r,"s"] & \mathcal{D} 
\end{tikzcd}
\]
In this diagram, $\varphi(z)=\rho z^{-1}$ and $\theta(z)$ is the conformal 
mapping 
which maps the set $\mathcal{G}^c$ onto the open disk
$\mathcal{D}$~\cite{KT-1985}:
\begin{equation}
  \label{eq:19}
  \theta(z)=\left(1-\sqrt{\frac{1-\frac{L-m}{2L}z}{1-\frac{m-L}{2m}z}}\right)
\left(1+\sqrt{\frac{1-\frac{L-m}{2L}z}{1-\frac{m-L}{2m}z}}\right)^{-1}.
\end{equation}

According to this diagram, $S(z)=\theta^{-1}(s(\rho z^{-1}))$,
$s(z)=\theta(S(\rho z^{-1}))$. The interpolation data for $s(z)$ are
obtained from~(\ref{eq:18}), 
$  s(0)=\theta(1)=\rho^*$, $ s(\rho)=\theta(0)=0$. 
A function $s(z)$ which solves the Nevanlinna-Pick problem with these data
exists if and only if 
\begin{equation}
  \label{eq:21}
  \left[
    \begin{array}{cc}
      1-(\theta(1))^2 & 1 \\
      1 &  \frac{1}{1-\rho^2}
    \end{array}
  \right]> 0;
\end{equation}
e.g., see~\cite[Theorem~18.1]{BGR-2013}. Condition~(\ref{eq:21}) is
equivalent to $\rho> \theta(1)=\rho^*$. This concludes the proof.
\hfill$\Box$


\begin{thebibliography}{10}

\bibitem{AM90}
B.~D.~O. Anderson and J.~B. Moore.
\newblock {\em Optimal Control: Linear Quadratic Methods}.
\newblock Prentice-Hall, 1990.

\bibitem{BGR-2013}
J.~Ball, I.~Gohberg, and L.~Rodman.
\newblock {\em Interpolation of rational matrix functions}.
\newblock Birkh{\"a}user, Basel, 1990.

\bibitem{dAST-2021}
A.~d'Aspremont, D.~Scieur, and A.~Taylor.
\newblock Acceleration methods.
\newblock {\em Foundations and Trends{\textregistered} in Optimization},
  5(1-2):1--245, 2021.
\newblock arXiv:2101.09545.

\bibitem{DT-2022}
Y.~Drori and A.~Taylor.
\newblock On the oracle complexity of smooth strongly convex minimization.
\newblock {\em Journal of Complexity}, 68:101590, 2022.

\bibitem{HC-2008}
W.~M. Haddad and V.~Chellaboina.
\newblock {\em Nonlinear dynamical systems and control: a {L}yapunov-based
  approach}.
\newblock Princeton University Press, 2008.

\bibitem{Khalil}
H.~K. Khalil.
\newblock {\em Nonlinear systems}.
\newblock Prentice Hall, Upper Saddle River, third edition, 2002.

\bibitem{KT-1985}
P.~Khargonekar and A.~Tannenbaum.
\newblock Non-euclidian metrics and the robust stabilization of systems with
  parameter uncertainty.
\newblock {\em IEEE Transactions on Automatic Control}, 30(10):1005--1013,
  1985.

\bibitem{Nesterov-2018}
Y.~Nesterov.
\newblock {\em Lectures on convex optimization}.
\newblock Springer, 2018.

\bibitem{OR-2000}
J.~M. Ortega and W.~C. Rheinboldt.
\newblock {\em Iterative solution of nonlinear equations in several variables}.
\newblock SIAM, 2000.

\bibitem{Polyak-1964a}
B.~T. Polyak.
\newblock Some methods of speeding up the convergence of iteration methods.
\newblock {\em USSR Computational Mathematics and Mathematical Physics}, 4(5):1
  -- 17, 1964.

\bibitem{Polyak-1987}
B.~T. Polyak.
\newblock {\em Introduction to optimization}.
\newblock Optimization Software, Inc., NY, 1987.

\bibitem{VFL-2017}
B.~Van~Scoy, R.~A. Freeman, and K.~M. Lynch.
\newblock The fastest known globally convergent first-order method for
  minimizing strongly convex functions.
\newblock {\em IEEE Control Systems Letters}, 2(1):49--54, 2017.

\end{thebibliography}
\end{document}